\documentclass[11pt,notitlepage,twoside,a4paper]{amsart}
\usepackage{amsmath}
\usepackage{amsfonts}
\usepackage{amssymb}
\usepackage{amscd}
\usepackage{amsthm}
\usepackage{amsbsy}
\usepackage[french, british]{babel}

\usepackage{amsmath,amssymb,enumerate}

\usepackage{epsfig,fancyhdr,color}%,showkeys,amsmidx?Aris

\usepackage{amssymb}
\usepackage{amsmath,amsthm}
\usepackage{latexsym}
\usepackage{amscd}
\usepackage{psfrag}
\usepackage{graphicx}
\usepackage[latin1]{inputenc}
\usepackage[all]{xy}
\usepackage[mathcal]{eucal}

\definecolor{NoteColor}{rgb}{1,0,0}

\renewcommand{\textsc}{\textcolor{red}}

\newtheorem*{theorem 1}{\rm\bf Proposition 1}
\newtheorem*{theorem 2}{\rm\bf Proposition 2}

\theoremstyle{definition}

\theoremstyle{remark}

\def\interieur#1{\mathord{\mathop{\kern 0pt #1}\limits^\circ}}

\begin{document}

 \begin{abstract}
 This is a biography and a report on the work of Vladimir Turaev. Using fundamental  techniques that are rooted in classical topology, Turaev introduced new ideas and tools that transformed the field of knots and links and invariants of 3-manifolds. He is one of the main founders of the new topic called quantum topology.  In surveying Turaev's work, this article will give at the same time an overview of an important part of the intense activity in low-dimensional topology that took place over the last 45 years, with its connections with mathematical physics. 
 The final version of this article appears in the book \emph{Topology and Geometry
A Collection of Essays Dedicated to Vladimir G. Turaev},   EMS Press, Berlin, 2021.

\bigskip

\noindent Keywords:  Vladimir Turaev ; Higher linking numbers;  intersection of loops on surfaces;
 Turaev cobracket;
Poincar\'e complexes; Poincar\'e duality;
 spin structures in 3-manifolds;
explicit constructions of cocycles;
 Turaev surface; Turaev volume;
skein module;
quantum invariants of links and 3-manifolds;
   $6j$-symbols; Turaev--Viro invariant;
Reshetikhin--Turaev invariant;  TQFT;
HQFT;
Gauss words and links;
enumeration problems in topology and group theory;
generalizations of the Thurston norm;
knotoids; knots; links; braids;
intersections and self-intersections  of loops on surfaces;
cobrackets;
combinatorial group theory; metric geometry; phylogenetics.

\bigskip


 \end{abstract}
%------------------------------------------------------------------------------------------------------------------
\title{Vladimir Turaev, friend and colleague}
\author{Athanase Papadopoulos}
%------------------------------------------------------------------------------------------------------------------

%\author{Athanase Papadopoulos}
%\address{Athanase Papadopoulos,  Universit{\'e} de Strasbourg and CNRS,
%7 rue Ren\'e Descartes,
% 67084 Strasbourg Cedex, France}
%\email{papadop@math.unistra.fr}
%

\date{}
\maketitle

%------------------------------------------------------------------------------------------------------------------
%---------------------------------

\tableofcontents
 
% 
% \section{Introduction}
% 
% 
% 

  \section{Introduction}
  
  Vladimir Turaev just turned 66.
  This article is a survey of part of his contribution to mathematics, and of the impact of his ideas. I have also included a short biographical sketch, in two parts, concerning his life in the Soviet Union, and his life in the West. I constructed the biography using  facts I learned from him during a 30-year old friendship.
Regarding the mathematical part, I am aware of the fact that my exposition will be at some places too short for readers who do not know much about Turaev's work, and it will be redundant for the experts.  I apologize in advance for both categories of mathematicians.  At the same time, this article will give an overview of an important part of the intense activity in low-dimensional topology that took place over the last 45 years, namely, invariants of knots, links and 3-manifolds, influenced by ideas from mathematical physics. Finally, I think that going through these developments is an illustration of how  the work of a single man can spread and transform the activity of a whole community.

  \section{Biography: Leningrad}
  Vladimir Georgievitch Turaev was born in Leningrad (Soviet Union), in 1954.
  His mother was a theater critic, and his father a producer at a puppet theater. He entered elementary school in 1961 and finished high school in 1970. The high school he attended was known under the name Boarding School No. 45. It was famous in Leningrad and was intended for gifted students. It was founded by the academicians M. V. Keldysh, I. G. Petrovsky and I. K. Koikin and was attached to Leningrad University. It included classes specialized in physics-mathematics or in chemistry-biology.

 After high school, Turaev enrolled in the Faculty of Mathematics and Mechanics of Leningrad State University (mathematics section), and he graduated from there in 1975. In 1976, he was hired as a researcher at the Steklov Mathematical Institute of the Academy of Sciences of the USSR, Leningrad Branch (LOMI) and he worked there until 1990, first with the status of a trainee (1976--1977), then as a junior researcher (1978--1984) and finally, from 1985 on, as a senior researcher. He obtained his ``kandidat" diploma (the Soviet analogue of the PhD) in 1979, from the Steklov Institute of Moscow. His scientific advisors were O. Ya. Viro and V. A. Rokhlin (the latter was Viro's advisor). The title of Turaev's thesis was ``Alexander--Fox invariants of 3-dimensional manifolds and Reidemeister torsions". In this thesis, he introduced a new torsion invariant which is more powerful than the Milnor torsion and which he called several years later (namely, in his book \emph{Torsions of 3-dimensional manifolds} published in 2002) ``maximal abelian torsion". In 1988, Turaev obtained the diploma of ``Doctor in physical and mathematical sciences", the Russian equivalent of the French Doctorat d'\'Etat (or the actual French Habilitation), again from the Steklov Mathematical Institute of Moscow.  The title of his dissertation was ``Classification problems in 3-dimensional topology". 
  
  From 1973 to 1975, while he was a university student, Turaev taught geometry at his former high school. Between 1979 and 1983, while working at LOMI, he taught topology at Leningrad's pedagogical institute. During these first years as a mathematician, he took part in several topology conferences, among which three international conferences: Moscow (1977), Leningrad (1982) and Baku (1987). 
  
   Turaev wrote his first paper in 1974, at the age of 20.  The theme was the relation between the Alexander polynomial and Reidemeister torsion of 3-dimensional manifolds. The paper was based upon ideas of Milnor. Turaev's interest in the question was stimulated by a question of Viro, which was motivated by the following fact: The classification of lens spaces of arbitrary dimension was completed in 1935 by K. Reidemeister in dimension 3, using (what, later on, became known as) Reidemeister torsion \cite{R1}. In the same year, W. Franz generalized the definition of torsion and carried out the classification in all dimensions (see \cite{Franz}). But in dimension 3, another proof of the same  classification was given by E. J. Brody in 1960, using the Alexander polynomial \cite{Brody}. The question Viro asked was to understand the relation between the two approaches. Turaev's paper appeared in print two years later, under the name \emph{Reidemeister torsion and the
Alexander polynomial} \cite{Reidemeister1976}.
   
  Starting from 1973, Turaev became actively involved in the remarkable mathematical bubbling that took place in Leningrad around the seminar organized by Rokhlin.  The topics discussed at that seminar included combinatorial topology, algebraic topology, differential topology, ergodic theory and real algebraic geometry. Several outstanding mathematicians were trained at what became known as the ``Rokhlin School". They include Y. M. Eliashberg, M. L. Gromov, A. M. Vershik, O. Ya. Viro, N. V. Ivanov and others.  Between the years 1973 and 1981 (the year Rokhlin retired), Turaev gave more than 35 lectures at Rokhlin's seminar. He also gave talks at the seminars of S. P. Novikov, V. I Arnold, and M. M. Postnikov in Moscow.  
   
 Between 1980 and 1990, Turaev was the scientific secretary of the mathematics section of the Leningrad club of scientists. The section's president was Vershik. Turaev's duty was mainly the organization of the monthly lectures for the members of the club and for those of the Leningrad Mathematical Society. In 1983, Turaev gave a short talk at the Warsaw ICM.

 In the year 1985, Perestro\"\i ka started in the Soviet Union. The same year, Turaev was allowed to travel to the West and he participated  in a conference in Oberwolfach. In the period that followed, he visited the University of Geneva (2 months in 1986), the University of Sussex (one month in 1987) and East Berlin (one week at the Mathematics Institute in 1988). In 1988, he stayed for one month at the University of Budapest. In 1989, he made an extended visit to Western Europe: one month at the University of Paris-Sud (Orsay), 2 weeks at the University of Marseille, one month at the University of Strasbourg and 2 months at the University of Bochum. During these stays, he gave talks in Toulouse, Grenoble, Lyon, Mannheim, Bonn, Heidelberg, G\"ottingen, Frankfurt and Geneva. In 1990, he spent 3 months at Ohio State University.  During this stay in the US, he lectured at Chicago, Stony Brook, Yale, Harvard, Brandeis, San Diego, MIT and Berkeley. In the summer of the same year, he participated in a conference on knot theory in Osaka, shortly before going to the ICM in Kyoto, where he was an invited speaker.

  The Kyoto ICM was the culmination of Turaev's Leningrad period. At the same time, it marked the end of that period, since after the congress, and after spending only a few days in Leningrad, Turaev left the Soviet Union and moved to France. 
  
   Among the four recipients of the Fields medal at the Kyoto ICM were Drinfeld, Jones and Witten, three names with whom Turaev's work is closely associated.   Presenting Witten's work at the Fields medalists' ceremony, Atiyah declared \cite{Atiyah-Work}: 
  \begin{quote}\small
  Witten's approach immediately shows how to extend the Jones theory
from knots in the 3-sphere to knots in arbitrary 3-manifolds. This generalization
(which includes as a specially interesting case the empty knot) had previously
eluded all other efforts, and Witten's formulas have now been taken as a basis
for a rigorous algorithmic definition, on general 3-manifolds, by Reshetikhin and
Turaev.
  \end{quote}
  
  Turaev's talk at that ICM was titled \emph{State sum models in low-dimensional topology}. At the congress and in the paper published in the proceedings \cite{State1991}, he described his approach to the construction of a 3-dimensional Topological Quantum Field Theory (TQFT). His method was based on  Drinfeld's theory of quantum groups.  At the beginning of his paper  \cite{State1991}, he declares  that the psychological barriers between physics and topology were aborted by the introduction in the latter of fundamental ideas of statistical mechanics by Kauffman and of quantum
physics by Witten. The works of Kauffman and Witten were motivated by Jones' discovery of a new polynomial invariant of links in the 3-sphere using von Neumann algebras. In the same introduction, Turaev points out the dichotomy between the Bourbaki style of writing mathematics and the new ``non-rigorous" methods inspired from mathematical physics. The latter methods were destined to have an enormous impact in the field of topology, and Turaev was one of the pioneers in using them to extract rigorous constructions.

 The late 1980s, that is, the last few years of the period in which Turaev worked at the Leningrad Steklov Institute, were for him the epoch of two decisive collaborations: with N. Yu. Reshetikhin and with O. Ya. Viro. Both collaborations had tremendous impacts on his later work and on mathematics in general. We shall review this later in this article.

 Shortly after Turaev left the Soviet Union, the name Leningrad disappeared: in 1991, after a referendum, the city recovered its original name, Saint Petersburg. The high school where Turaev studied was renamed ``Academic Gymnasium" as it became a department of Saint Petersburg State University.
 In the same year, the Soviet Union was dissolved, and the Russian Federation became its successor state.  
 
  Before passing to the rest of the biography, I will make a short review of Turaev's mathematical work during the Leningrad period of his life.

    \section{Mathematics: 1973--1990}
    
   In the period 1973--1990, the research activity of Turaev was centered on the study of loops on surfaces, knot theory and 3-manifolds, three themes that constitute until now his favorite field of investigation. The ideas he presented in the papers he wrote during these years contain the germ of the work he developed in the years that followed.
   
   More specifically, during his Leningrad period, Turaev concentrated on the following 10 topics, roughly classified in chronological order:
   \begin{enumerate}
   \item  Reidemeister torsion;
      \item Higher linking numbers for links in the 3-sphere;
      \item Intersection of loops on surfaces;
      \item The Turaev cobracket;
         \item Poincar\'e complexes, Poincar\'e duality and its influence on the fundamental group;
            \item Spin structures in 3-manifolds;
   \item Explicit constructions of cocycles;
   \item Turaev surfaces and Turaev volume;
      \item Skein modules;
   \item Quantum invariants of links and 3-manifolds.
   \end{enumerate}

   Let me make a few comments on each of these topics. I will mention a few results of Turaev for each of them, leaving aside several others which perhaps are not less  interesting.
   
  \subsection{Reidemeister torsion} In his work on this subject, Turaev was motivated by a theorem of Milnor that appeared in his 1962 paper \emph{A duality theorem for Reidemeister torsion} \cite{Milnor1962}, stating that the Alexander polynomial of a link in a closed 3-manifold can be interpreted as a kind of Reidemeister torsion. Milnor used this result to obtain properties of the Alexander polynomial from universal properties of Reidemeister torsion. Turaev developed this theory in several directions. 
In his paper \emph{Towards the topological classification of geometric 3-manifolds}, published in 1988 \cite{Towards1988}, he defined an invariant for a class of oriented 3-manifolds which (since Thurston's geometrization conjecture became a theorem) we know is the class of all oriented closed 3-manifolds. 

In another paper, titled \emph{Euler structures, nonsingular vector fields, and Reidemeister-type torsions}, published in 1989  \cite{Euler1989}, Turaev applied this theory to the study of the action of diffeomorphisms of smooth manifolds on nonsingular vector fields. In the same paper, he introduced the notion of \emph{Euler structure},  as an equivalence class of certain non-singular vector fields on a smooth manifold which may also be defined using some special chains (called \emph{Euler chains}) on a PL triangulation of the manifold. The concept of Euler structure  leads to a notion of torsion that refines the Reidemeister torsion. In fact, together with the notion of homology orientation, an Euler structure controls the indeterminacy of the Reidemeister torsion. 

Several years later, Turaev used Euler structures in a paper he wrote with M. Farber titled \emph{Poincar\'e--Reidemeister metric, Euler structures, and torsion}, published in 2000 \cite{Poincare2000}. Based on this notion, Farber and Turaev were able to introduce the so-called  ``Poincar\'e--Reidemeister scalar product''  on the determinant line of the cohomology of any flat vector bundle over a closed orientable
odd-dimensional manifold.

 \subsection{Higher linking numbers for links in the 3-sphere}  This work also originates in an idea of Milnor.  The latter introduced these higher linking numbers in 1957, in his paper \emph{Isotopy of links} \cite{Milnor1957}, as a generalization of Gauss's linking coefficients for closed curves in 3-space. Among Turaev's results on this topic, we mention the paper \emph{The Milnor invariants and Massey products} \cite{Milnor1976} (1976) in which he gave a cohomological description of higher linking invariants using Massey products in link complements.
 
 Together with Milnor's name, that of Gauss is important here. The latter's work on links acted later on as an inspiration for a series of works of Turaev. We shall say more about this in \S \ref{s:Gauss} below. 

   \subsection{Intersection of loops on surfaces} 
Among Turaev's works on this subject, we mention the paper \emph{Intersections of loops in two-dimensional manifolds} 
\cite{Intersections1978} published in 1978,  in which he defined a homotopy intersection form for curves on surfaces. In this article, he obtained several results, including necessary and sufficient conditions for an element of the fundamental group of a surface to be realizable by an embedding of a circle.  This is a 2-dimensional analogue of a higher-dimensional result of Kervaire.  The introduction of the homotopy intersection form was motivated by its applications in knot theory, where analogous intersection forms are defined on knot groups. At the same time, Turaev noticed that a version of this form is inherent in the work of Papakyriakopoulos (\emph{Planar regular coverings of orientable closed surfaces} \cite{Planar1975}, 1975).  Incidentally, in the same paper, Turaev gave a new proof of the fact that an automorphism of the fundamental group of a surface with nonempty boundary preserving the peripheral structure is realized by a self-homeomorphism of the surface (this is usually called the Dehn--Nielsen--Bauer theorem).  

A sequel to this paper, carrying almost the same title 
    (\emph{Intersection of loops in two-dimensional manifolds. II. Free loops})  appeared  in 1983 and was written jointly with Viro \cite{Intersections1983}. The paper is concerned with the problem of finding the minimal  number  of intersections and self-intersections of a loop in a given homotopy class on a surface.
  
   Several years later, Turaev's form played an important role in the theory of representations of fundamental groups of surfaces; see e.g. the works by  G. Massuyeau and Turaev \cite{Fox2013} (2013) and \cite{Quasi2014} (2014), and \S \ref{s:2-intersection} below.

  \subsection{The Turaev cobracket}   \label{cobracket1}
In the paper \emph{Algebras of loops on surfaces, algebras of knots, and quantization} \cite{Algebras1989}, Turaev introduced a Lie cobracket on the free module generated by the homotopy classes of loops on an oriented surface. Moreover, he noticed that together with the Lie bracket introduced by W. Goldman in his paper \emph{The Symplectic Nature of Fundamental Groups of Surfaces} \cite{Goldman1984} (1984), the cobracket turns this free module into a Lie bialgebra in the sense of Drinfeld \cite{Drinfeld1983}. Turaev constructed a biquantization of this Lie bialgebra.
   He further developed this theory in his paper titled \emph{Skein quantization of Poisson algebras of loops on surfaces} \cite{Skein1991},  published in 1991, in which the notions of bi-Poisson bialgebra and biquantization took shape.
   
    The question of the geometric interpretation of the Turaev cobracket arose naturally, especially since it was known that the Goldman bracket has such an interpretation: it is the Poisson bracket defined by the Weil--Petersson symplectic structure of Teichm\"uller space. Turaev's cobracket had several developments.     In a paper titled \emph{Loops on surfaces, Feynman diagrams, and trees} \cite{Loops2005}, published in 2005, Turaev gave a geometric meaning of his cobracket, making a relation with the Lie cobracket defined by Connes and Kreimer in their work on perturbative quantum field theory \cite{Connes-Kreimer1998}. 
    
    Turaev's cobracket played a major role in a series of papers by him and  Massuyeau published between 2013 and 2018 \cite{Fox2013, Quasi2014, Brackets2017, Brackets2018}, and  in other papers published during the same period by N. Kawazumi, Y. Kuno, A. Alekseev and F. Naef \cite{Higher2017, Intersection2015, Goldman-Turaev2018}. We shall mention some of these developments in  \S \ref{cobracket2} below. A detailed exposition of the Turaev cobracket together with some variants and applications is contained in the chapter by Kawazumi in the present volume \cite{Kawazumi-Turaev}.

   \subsection{Poincar\'e complexes, Poincar\'e duality and its influence on the fundamental group} Poincar\'e complexes are CW-complexes that are homotopy analogues of closed manifolds. Their study was initiated by S. P. Novikov and W. Browder. A question that arose since the birth of this theory was the  characterization of Poincar\'e complexes that are homotopy-equivalent to 3--dimensional closed manifolds. Turaev worked on this problem and obtained  a series of results that led to a complete algebraic reduction of the homotopy classification of 3-dimensional Poincar\'e complexes. Among the results of his paper  \emph{Three-dimensional Poincar\'e complexes: homotopy classification and splitting} \cite{Three-dimensional1989} published in 1989, there is an answer to a question posed by C. T. C. Wall in 1967, asking for a generalization to 3-dimensional Poincar\'e complexes of a classical splitting theorem of Kneser and Stallings for 3-manifolds, see \cite{Wall1965}.

 \subsection{Spin structures in 3-manifolds}  This is one of the topics on which Turaev worked which is profoundly rooted in the Rokhlin school. Spin structures on a smooth manifold are elements of the first cohomology group with coefficients in $\mathbb{Z}/2$ of the manifold of positively oriented bases in the tangent bundle of the given manifold. (The cohomology class is assumed to be non-trivial on the fibers.) Every oriented 3-manifold $M$ admits exactly $2^s$ spin structures, where $s= \mathrm{dim} \ H_1(M;\mathbb{Z}/2)$.  (In fact, there is a natural affine parametrization of spin structures by $H_1(M;\mathbb{Z}/2)$.) The definition of the famous Rokhlin invariant is based on the fact that the signature of the intersection form of a spin smooth closed 4-manifold is divisible by 16 (this is Rokhlin's theorem).

 Turaev  used spin structures in the classification of isotopy types of oriented links in the 3-sphere. In his paper \emph{Classification of oriented Montesinos links via spin structures}  \cite{Classification1988} (1988), he classified oriented Montesinos links after establishing a canonical correspondence between the orientations of a link in the 3-sphere, up to reversal of all the orientations, and the spin structures on its two-fold branched cover. Among his other publications on this topic, we mention the paper \emph{Cohomology rings, linking coefficient forms and invariants of spin structures in 3-dimensional manifolds} \cite{Cohomology1983} (1983) in which he provided necessary and sufficient conditions on a triple consisting of a sequence of graded rings, a bilinear form, and a function with values in $\mathbb{Z}/16$ to be respectively the cohomology ring, the linking form, and the Rokhlin function of some closed oriented 3-manifold.  The Rokhlin function  of a manifold is precisely the function which assigns to a spin structure on this manifold its Rokhlin invariant.

 \subsection{Explicit construction of cocycles}  In his two papers \emph{A cocycle of the symplectic first Chern class and Maslov indices}   \cite{Cocycle1984} (1984) and \emph{The first symplectic Chern class and Maslov indices} \cite{First1987} (1987), Turaev  introduced a 2-dimensional Borel cocycle on the symplectic group  $\mathrm{Sp}(n)$ and established a relation between this cocycle and the Maslov indices of linear Lagrangian spaces. The result sheds a new light on important formulae of Leray concerning the Maslov indices of the Maslov class and it gives a novel way of computing them.  Several years later, Endo and Nagami made a relation between the cocycle introduced by Turaev and the so-called Meyer cocycle \cite{EN}.
 
 \subsection{Turaev surfaces and Turaev volume} 

In a paper published in 1987, titled  \emph{A simple proof of the Murasugi and Kauffman theorems on alternating links} \cite{Simple1987},  Turaev gave a new proof of the fact that the span of the Jones polynomial
of a nonsplit alternating link is equal to its crossing number. In doing so, he introduced a notion that was revived several years later under the name \emph{Turaev surface} and which plays today a significant role in low-dimensional topology.  It leads to an invariant, called the \emph{Turaev genus} of a knot, obtained by taking the minimal genus of a Turaev surface of this knot.  Roughly speaking, Turaev surfaces play for the Jones polynomial the role which is played by   Seifert surfaces for the Alexander polynomial. These surfaces give a measure of how far a knot is from being alternating. They also appear  in other settings; see for instance the paper \cite{DL} in which the authors develop a Turaev surface approach to Khovanov homology.

There is an expository article on  Turaev surfaces in the \emph{Concise Encyclopedia of Knot Theory} \cite{KK}.   
   
In the recent paper \emph{Turaev Hyperbolicity of Classical and Virtual Knots} by 
 Adams, Eisenberg, Greenberg, Kapoor, Liang, O'Connor, Pacheco-Tallaj and  Wang   \cite{Hyperbolicity} (2019), the authors introduce the notion of \emph{Turaev volume} of an arbitrary virtual or classical link, as the minimal volume among all the hyperbolic 3-manifolds associated via  Turaev's construction which originates in his paper \cite{Simple1987}.

\subsection{Skein modules} A skein relation is a local relation involving several different ways of modifying a link inside a small ball. Such relations are used to obtain polynomials associated to knots. Skein modules and algebras are defined using such polynomials. The theory of skein modules of an oriented 3-manifold was developed by Turaev,  who called them Conway and Kauffman modules, in his paper \emph{The Conway and Kauffman modules of a solid torus} \cite{Turaev1988} (1988). In this paper, the Conway module of an oriented 3-manifold  is defined as the module of formal linear combinations of isotopy classes of oriented links in this manifold with coefficients in the ring of Laurent polynomials in 2 variables, $\Lambda=\mathbb{Z}[x,x^{-1}, y, y^{-1}]$, considered up to local relations known as Conway's relations. The Kauffman module is defined similarly for isotopy classes of framed links, with the local relations being Kauffman's relations. In his paper   \cite{Turaev1988} , Turaev calculates the Conway module of a solid torus $\mathbb{S}^1\times I^2$.  In this work, Turaev uses a Yang--Baxter equation and the $R$-matrices which arise in this theory. Josef Przytycki gave an independent development of an analogous theory in his paper \emph{Skein modules of 3-manifolds} \cite{Przytycki1991} (1991).

 \subsection{Quantum invariants}  \label{quantum1} The first paper by Turaev on this topic, \emph{The Yang--Baxter equation and invariants of links} \cite{Yang1988}, was published in 1988.  It was the starting point of a series of publications by him  that had tremendous consequences in the three decades that followed. In this paper, Turaev developed a general scheme which produces isotopy invariants of links in the 3-sphere, using $R$-matrices from statistical mechanics, i.e., solutions of the fundamental Yang--Baxter equation. He showed that this scheme includes the Jones polynomial and its generalizations: the Kauffman and the HOMFLY (now called HOMFLY-PT) polynomials.  Ideas from statistical mechanics applied to the Jones polynomial were already brought into the theory of the Jones polynomial by Jones himself and by Kauffman, and 
Turaev's paper made it clear that these generalizations of the Jones polynomial are also closely related to vertex models in statistical mechanics.    In the same paper, Turaev developed a general theory of tangles, that is, embeddings of a certain number of arcs and circles in the 3-ball. The theory of tangles allows operations on tangle diagrams which lead to categories and algebraic structures that turned out to be among the main tools in the study of knots and links.

In a paper published in the following year, \emph{Operator invariants of tangles and $R$-matrices} \cite{Operator invariants1989}, Turaev introduced a new category of tangles that generalize both representations of braid groups involving $R$-matrices and the Jones--Conway and Kauffman polynomials of links that were recently introduced. Turaev's construction is based on the notion of quantum $R$-matrix.

After he wrote these two papers, Turaev started his collaboration with Reshetikhin, with whom he generalized this theory to colored tangles. 
The latter had just written a paper on the same subject titled \emph{Quantized universal enveloping algebras, the Yang--Baxter equation, and invariants of links} \cite{R1988} (1988). Turaev and Reshetikhin's joint paper \emph{Ribbon graphs and their invariants derived from quantum groups} \cite{Ribbon1990} appeared in 1990. Here, a ``color" for a tangle is an $A$-module associated to it, where $A$ refers to a quasitriangular Hopf algebra. The quantum invariant appears  as a functor from a category of diagrams where the tangles act as morphisms,  to the category of A-modules.

 In a subsequent paper, \emph{Invariants of 3-manifolds via link polynomials and quantum groups} \cite{Invariants1991}, published in 1991, Reshetikhin and Turaev constructed a 3-dimensional TQFT that provides nontrivial invariants of 3-manifolds from a modular category, introducing in this context  the notion of modular Hopf algebra (a quantum group), using the invariants of framed links introduced in their previous paper 
 \cite{Ribbon1990} and surgery techniques based on Kirby's theorem \cite{Kirby}. With this work, the Jones polynomial, initially defined for links in the 3-sphere, was generalized to links in arbitrary compact oriented 3-manifolds. At the same time, the paper by Reshetikhin and Turaev makes an important step in the program initiated by Witten that produces invariants of 3-manifolds using quantum and conformal field theories. Before that, Witten had given an interpretation of the Jones polynomial in terms of quantum field theory \cite{Witten1989}, and the work done by Turaev and Reshetikhin gave a mathematical model for this theory, based on Drinfeld's theory of quantum groups and using the classical Kirby calculus of links in the 3-sphere. An overview of some of the important ideas in this work is contained in the paper \cite{BDR} by Blanchet and De Renzi in the present volume.

At the same time, Turaev started a new collaboration with Viro on the construction of quantum-type invariants of 3-manifolds based on quantum $6j$-symbols  and state sum models of statistical mechanics (analogues of partition functions). This led to the paper 
 \emph{State sum invariants of 3-manifolds and quantum $6j$-symbols} \cite{State1992}, published in 1992, where a new invariant for 3-manifolds is defined as a certain state sum computed using an arbitrary triangulation of the manifold. The construction involves  summation over certain colorings (maps from the  set of edges of the triangulation  to some finite set) and  the so-called $6j$-symbols arising in the representation theory of quantum groups. A substantial part of the construction of the invariant consists in showing that it is independent of the choice of the triangulation and the coloring. This is proved using combinatorial techniques developed by M. H. A Newman and J. W. Alexander in the 1920s. The ideas in the paper do not only lead to numerical invariants of manifolds, but also to a 3-dimensional TQFT. In other words, these invariants behave in a functorial manner with respect to the gluing of manifolds along their boundaries.
 
 In the paper \emph{Shadow links and face models of statistical mechanics} \cite{Shadow1992} which appeared in 
the same year, Turaev developed a new approach to quantum invariants of links in circle fiber spaces over surfaces. In doing so, he introduced a geometrical representation of links as ``shadows", that is, systems of loops on surfaces, with an integer attached to each component of the complement of a loop. This gave a new method of presenting links in the 3-sphere, using shadows on the 2-sphere. Here, the 3-sphere is considered as fibered by the Hopf fibration over the 2-sphere. The diagrammatic isotopy which lifts to an ambient isotopy of links is based on a set of moves that are analogues of the Reidemeister moves. Using complex numbers instead of integers, Turaev leads us to the notion of complex shadow links. 
Models of statistical mechanics were used again by him to construct non-trivial Jones-type invariants of shadow links. 

Soon afterwards, in a short note titled \emph{Quantum invariants of 3-manifolds and a glimpse of shadow topology} \cite{Quantum1991} (1991), Turaev explained the relation between the Reshetikhin--Turaev invariant defined by surgery and quantum groups, and the  Turaev--Viro invariant defined using state sums on triangulations. The same relation was independently discovered by K.  Walker \cite{Walker}. This relation was also touched upon in Turaev's ICM talk \cite{State1991}. A large section of his book \emph{Quantum invariants of knots and 3-manifolds} \cite{Quantum1994} which appeared in 1994 is dedicated to this question. Turaev's theory of shadows is also reviewed in this book

Several years later, Turaev's ideas on shadows were taken up by F. Costantino and D. Thurston in their paper \cite{CoT} in which they introduced a notion of \emph{shadow complexity} for compact orientable 3-manifolds with (possibly empty) toral boundary. They used this notion to prove several results concerning volumes of hyperbolic manifolds, including a double inequality relating this quantity and the Gromov norm.

  \section{Biography: The West}
 
  In April 1990, Turaev was offered a position of Directeur de
Recherche at the French CNRS (Centre National
de la Recherche Scientifique). This is a research position, at the senior level, free from any teaching, equivalent to the one he had at the Steklov Institute in Saint Petersburg. His career as a researcher in France started on the first of September of the same year, and he chose to settle in Strasbourg and to work at the Institut de Recherche Math\'ematique Avanc\'ee (IRMA) there. In the meanwhile, he got several attractive offers from American universities and he declined them all. His new research position in France allowed him to spend several months abroad each year. He gave lectures in various countries including Switzerland, Denmark, and the US. He also made several long stays in Japan. In Strasbourg, he directed several PhDs and he stimulated the activities of several researchers, young and older. The list of his new collaborators included C. Kassel, G. Massuyeau, A. Virelizier, N. Geer and several others.

   In 2004, Vladimir Turaev was awarded the silver medal of CNRS. This is the most prestigious prize given to researchers in their mid-career in France.
 Between the years 2000 and 2008, he  was in charge of the ``Rencontres entre math\'ematiciens et physiciens th\'eoriciens", a bi-annual 3-days meeting in Strasbourg, one of the oldest mathematical regular meetings worldwide which is still operating. In 2000, he founded the IRMA Lectures in Mathematics and Theoretical Physics, a collection of books first published by De Gruyter and later on by the European Mathematical Society, of which he was the editor-in-chief, until 2019. In fact, it was an idea that I heard of him since he arrived to France that the community of mathematicians is in need of good survey papers, and several volumes in the IRMA series consist of collections of survey papers. Turaev served also on the editorial board of several journals,  and I would like to mention especially the journal Quantum Topology which he founded in 2009 and of which he is editor-in-chief.  The journal is published by EMS; the first issue appeared in 2010.

In 2008, Turaev was appointed   at the University of Indiana where he has  taught as the Boucher Professor, while keeping  close relations with Strasbourg. 

 In 2014, Turaev returned to Russia on a part-time basis, as the director of the Laboratory of Quantum Topology  based at Chelyabinsk State University in Siberia. This laboratory included groups of scientists working in Chelyabinsk, Moscow, Novosibirsk and Saint Petersburg. The project was funded by a grant of the Government of the Russian Federation for a period of 4  years.

        \section{Mathematics: 1990--today}
  
  After 1990, Turaev continued developing the topics which he used to work on in Leningrad. At the same time, new subjects of interest arose, novel ideas came to his mind, and they led to beautiful chapters in mathematics. His work  stimulated a large number of researchers, he certainly was also stimulated by some of them, and he had several new co-authors.
  I will review some of the topics on which he worked since 1990.
  I have classified them in the following 11 categories:
  \begin{enumerate}
    \item  $6j$-symbols, Turaev--Viro invariants and applications;
     \item  Reshetikhin--Turaev invariants and generalizations;
  \item  TQFT and unoriented TQFT;
    \item  HQFT;
  \item  Gauss words and links;
  \item Enumeration problems in topology and group theory;
   \item  Generalizations of the Thurston norm;
  \item  Knotoids, knots, links and braids;
   \item Intersections and self-intersections  of loops on surfaces;
  \item   Cobrackets;
  \item   Combinatorial group theory, metric geometry and phylogenetics.
  \end{enumerate}

\subsection{$6j$-symbols, Turaev--Viro invariants and applications}
 
Turaev published several papers in which he extended his work with Viro which resulted in the Turaev--Viro invariant based on $6j$-symbols.  In his paper \emph{Quantum invariants of links and 3-valent graphs in 3-manifolds} \cite{Quantum1993} (1993),  he introduced the notion of  quasimodular Hopf algebra, giving a new interpretation of the Jones polynomial of framed links and its generalizations. This interpretation is 3-dimensional, that is, it is not based on a 2-dimensional link projection, but on a study of the  link complement, using a refined version of the Turaev--Viro state sums on triangulations and $6j$-symbols. He also gave a 3-dimensional way of computing the values of this polynomial at roots of unity in terms of the link complement. A 3-dimensional interpretation of the Jones polynomial was already given by Witten, using Feynman integrals. Turaev writes in the introduction of this paper: ``Though the main
ideas of the paper are independent of Witten's approach one may view this paper as
an attempt to understand his work." Compared to Witten's approach, Turaev-Viro's theory, based on $6j$-symbols, is  mathematically rigorous. 
 The theory of $6j$-symbols also constitutes an important theme of Turaev's book \emph{Quantum invariants of knots and 3-manifolds} \cite{Quantum1994}.

Among the other works on the same subject, we mention Turaev's paper \emph{Tetrahedral forms in monoidal categories and 3-manifold invariants} \cite{Tetrahedral2012}  (2012) written in collaboration with N. Geer and R. Kashaev in which the authors produce new Turaev--Viro type invariants of 3-manifolds via monoidal categories. We also mention the  paper \emph{The Tambara-Yamagami categories and 3-manifold invariants} \cite{Tambara2012} (2012) by  Turaev and L. Vainerman in which these authors study generalized Turaev--Viro invariants of closed 3-manifolds arising from a tensor category. 

We also mention the two papers \emph{Modified quantum dimensions and re-normalized link invariants} \cite{Modified2009}  (2009) and \emph{Modified $6j$-symbols and 3-manifold invariants}  \cite{Modified2011} (2011) by Turaev, N. Geer and B. Patureau-Mirand in which these authors produce new 3-manifold state sum invariants using a  theory of modified $6j$-symbols

The paper  \emph{Three-dimensional manifolds with poor spines} \cite{T-Vesnin} by Turaev, E. Fominykh and A. Vesnin is concerned with problems of complexity of 3-manifolds in the sense of S. Matveev \cite{Matveev}. These authors use the notions of ``special polyhedron" and ``special spine" which originate in combinatorial topology, a terminology introduced by Matveev.  A special spine is said to be poor if it
does not contain proper simple sub-polyhedra. Using a particular case of the Turaev--Viro
invariants, called the $\epsilon$-invariant, they prove that if a compact 3-manifold
with connected nonempty boundary has a poor
special spine with 2 components and $n$ true vertices, then its complexity is equal to $n$. 
 They develop the theory of poor special spines of such a 3-manifolds and they construct examples of such 3-manifolds for an infinite number of values of $n$.

\subsection{Reshetikhin--Turaev invariants and generalizations}

Among the papers in which Turaev further generalized or improved the Reshetikhin--Turaev construction of 3-manifold invariants, we mention the article \emph{ Modular categories and 3-manifold invariants} \cite{Modular1992} (1992) in which he introduces in the representation theory of quantum groups the notion of modular ribbon category.  The notion of modular category is very close to that of modular tensor category in the sense of Moore and Seiberg. In his paper, Turaev shows that the notion of modular category suffices for the Reshetikhin-Turaev construction of TQFT in 2+1 dimensions, and therefore, it leads to invariants of 3-manifolds via surgery presentations. 

Among the papers concerned with generalized topological invariants of closed oriented 3-manifolds, we mention Turaev's paper with H. Wenzl titled \emph{Quantum invariants of 3-manifolds associated with classical simple Lie algebras}, published in 1993 \cite{Quantum1993-Wenzl}, in which the authors rely on the notion of quasimodular Hopf algebra to obtain 3-manifold invariants associated with classical simple Lie algebras.

\subsection{TQFT, unoriented TQFT  and HQFT}
Let us recall first that a TQFT is a theory which assigns, in a functorial way,  topological invariants to manifolds, using ideas from quantum field theory. More precisely, for $d\geq 1$, a $d$-dimensional TQFT over a commutative ring $K$ is an assignment to each closed oriented $d$-manifold $M$ a $K-$module $A_M$ and to every  compact oriented $(d+1)$-dimensional cobordism  $(W,M,M_1)$ a $K$-homomorphism $\tau(W):A_M\to A_{M_{1}}$. The construction must satisfy a certain set of axioms, which differ among authors, and this is why there are different TQFTs.
  After the notion of TQFT was introduced by Witten in the interpretation of the Jones polynomial, several systematic studies  of TQFTs and of the associated idea of modular functors were conducted by  Segal, Moore--Seiberg and Atiyah, who provided various sets of axioms for TQFTs.
  
   In his paper \emph{Axioms for topological quantum field theories} published in 1994 \cite{Axioms1994}, Turaev gave a new set of axioms for modular functors and TQFTs. In the introduction to his paper, he explains the setting: 
   \begin{quote}\small
   The problem with any axiomatic definition is that it should be sufficiently general but not too abstract. It is especially hard to find the balance in axiomatic systems for TQFTs because our stock of non-trivial examples is very limited. There is no doubt that further experiments with axioms for TQFTs will follow. The reader will notice that our definitions and results have a definite flavor of abstract nonsense. However, they form a natural background for more concrete 3-dimensional theories.
   \end{quote}
   
   In the same year, Turaev published his book  \emph{Quantum invariants of knots and 3-manifolds} \cite{Quantum1994}, in which he presented in detail his constructions with Reshetikhin and Viro of 3-dimensional TQFTs derived  from modular categories.  A review of this book, by G. Kuperberg, was published in the Bulletin of the AMS  \cite{Kuperberg}.

Turaev returned to this theme several years later.  In his paper \emph{Dijkgraaf--Witten invariants of surfaces and projective representations of groups} \cite{Dijkgraaf2007} (2007), he reformulated a TQFT theory proposed by Dijkgraaf and Witten in their paper 
  \cite{DW}  which is valid for $n$-dimensional manifolds, based on physical considerations and associated to compact and non-necessarily connected Lie groups $G$. He considered the case where the group $G$ is finite, and he gave a version of the Dijkgraaf--Witten invariant of surfaces in terms of projective representations of the group $G$. In a subsequent joint paper with S. Matveev \cite{Matveev1}, he expressed the Dijkgraaf--Witten invariants in terms of the Arf invariant, for the case where $G=\mathbb{Z}/2\mathbb{Z}$. King, Matveev, Tarkaev, and Turaev, in the more recent work \emph{Dijkgraaf--Witten  $\mathbb{Z}/2$-invariants for Seifert manifolds} \cite{King} (2017), gave precise formulae for this theory in the case of Seifert 3-manifold with orientable base.

In his paper \emph{Sections of fiber bundles over surfaces and TQFTs} \cite{Sections2010} (2010), Turaev used a TQFT in the solution of an enumeration problem for sections of Serre fibrations over compact orientable surfaces.

In their paper \emph{Unoriented topological quantum field theory and link homology} \cite{Unoriented2006} (2006), Turaev and Turner introduced the notion of unoriented TQFT and explained how it yields a link homology  theory   for link diagrams
on orientable surfaces which is invariant under stable equivalence classes of diagrams on surfaces, that is, under the equivalence relation generated by homeomorphisms of the surface, Reidemeister moves and the addition or subtraction of handles disjoint from the diagram. 

   \subsection{HQFT}
  In 1999, Turaev introduced the notion of Homotopy Quantum Field Theory (HQFT) as a version of TQFT  where the basic objects are manifolds endowed with homotopy classes of maps to a fixed topological space (called the ``target" of the theory). In his preprint  \emph{Homotopy field theory in dimension 2 and group-algebras} \cite{Homotopy1999} (1999),  starting with ideas from TQFT, he develops at the same time the general setting of HQFTs and  the algebraic structures underlying them. He introduces, for a given group $\pi$ a notion of $\pi$-algebra. He discusses lattice models for $(1+1)$-dimensional HQFTs with target $K(\pi, 1)$, developing the theory of cohomological  HQFTs with target $K(\pi,1)$ derived from cohomology classes of $\pi$ and its subgroups of finite index. He classifies (1+1)-dimensional HQFTs in terms of so-called crossed group-algebras and in particular,  he classifies the cohomological (1+1)-dimensional HQFTs over a field of characteristic 0 by simple crossed group-algebras.  At the same time, he introduces two state sum models for (1+1)-dimensional HQFTs and he proves that the resulting HQFTs are direct sums of rescaled cohomological HQFTs. In the same paper, he discusses deformations of Frobenius algebras and he deduces a version of the Verlinde formula from the fact that the lattice $(1+1)$-dimensional HQFTs over algebraically
closed fields of characteristic 0 are semi-cohomological. 
  
  One year after the preprint \cite{Homotopy1999}, Turaev wrote another paper on the same subject, \emph{Homotopy field theory in dimension 3 and crossed group-categories} \cite{Homotopy2000}, which sets the basis of the theory in dimension 3.
  
HQFT was further developed by Tim Porter in \cite{Porter2007}, by Porter and Turaev in \cite{PT2007}, by Staic and Turaev in \cite{ST2010} and by others. Turaev's monograph \emph{Homotopy quantum field theory}  \cite{Homotopy2010} published in 2010 supersedes the two preprints \cite{Homotopy1999} and  \cite{Homotopy2000} which we mentioned above. It sets in a more formal manner the foundations of this theory, using the notion of graded monoidal category and  involving several  new algebraic structures. 
   This monograph emphasizes again the surface case. A review of this book, by Porter, was published in 2012  in the Bulletin of the AMS \cite{Porter2012}.
   
    In a series of subsequent papers  with Virelizier \cite{TV1, TV3, TV2, TV4}, Turaev developed more fully the 3-dimensional case. This is the occasion for me to mention the book by Turaev and Virelizier, \emph{Monoidal categories and topological field theory} \cite{TV},  published in 2017,  in which the authors prove an important conjecture that makes a relation between the Turaev--Viro and the Turaev--Reshetikhin invariants. The book obtained in 2016 the Ferran Sunyer i Balaguer Prize for a mathematical monograph of an expository nature.

  \subsection{Gauss words and links} \label{s:Gauss}
  
    This section, central in the flow of the present paper, besides surveying Turaev's work on the theory of Gauss words, is meant to recall that knot and link theory, with its developments  based on the study of plane projections and diagrams of knots and links, is rooted in the work of Gauss.
    
A \emph{Gauss word}  is a sequence of letters in some given finite alphabet. The word is defined up to a circular permutation, in which all letters of the alphabet occur exactly
twice. A generic closed curve in the plane
gives rise to a Gauss word, obtained by representing this curve by a circle, labeling its self-crossings by different letters and writing them down
in the order of their appearance when the curve is  traversed. Gauss introduced this concept in his study of  the combinatorics  of closed curves in the plane, and he  gave a necessary condition for a Gauss word to be realized by a closed curve in
the plane  \cite[p. 272, 282--286]{Gauss}. Necessary and sufficient conditions for the same property to hold were obtained later by several authors. 

Based on these ideas of Gauss, Turaev developed a general topological theory of words
in his paper \emph{Topology of words} \cite{Topology2007}  (2007). In this work, general words are approximated by Gauss words and are studied up to a set of transformations that are 
inspired by homotopy of curves in the plane. The notion of homotopy of curves can  be translated into some
 local transformations of words  and  generate an equivalence relation called homotopy of words.
 Turaev  writes in the introduction: 
 \begin{quote}\small
 Words are finite sequences of letters in a given alphabet. Every word has
its own personality and should be treated with the same respect and attention
as say, a polyhedron or a manifold. In this paper we attempt to study words as
topological objects.
\end{quote}
In the same paper, he introduced the concepts of \'etale words and
nanowords. A letter appearing
in a nanoword occurs exactly twice.  Inspired by the usual methods of topology (self-linking function,
 linking form, linking pairing,  Alexander matrices,  modules and polynomial invariants, coloring, knot quandles, virtual strings, etc.),  Turaev developed several methods for the study of the properties of words that are preserved by homotopy.
 
In his paper \emph{Knots and words} \cite{Knots2006} (2006), Turaev realized a unification of several aspects of knot theory using Gauss words,  extending the Kauffman
bracket polynomial, the Jones polynomial, the knot quandle and other notions and methods to the more general setting of words and phrases in arbitrary alphabets. For this purpose, he introduced the notion of stable equivalence classes of
(pointed oriented) knot diagrams on a surface. These equivalence classes are in one-to-one correspondence with (appropriately defined)
homotopy classes of nanowords in an alphabet consisting of four letters.    
Classical and virtual links are interpreted as nanophrases.

Turaev's paper \emph{Virtual strings}  \cite{Virtual2004}  (2004) is also concerned with Gauss words. In this paper, the notion of virtual string is introduced. This is a scheme of self-intersection of an oriented generic closed curve on an oriented surface, represented by an oriented circle  equipped with a certain number of ordered pairs of points that are called arrows.  A  virtual string gives rise to a Gauss word  by labeling the arrows with different letters (the ``alphabet") then traversing 
the circle  in the positive direction and keeping track of  the label of arrows.   
This notion of virtual string is  inspired from the notion of virtual knot introduced by  Kauffman in his paper \emph{Virtual knots theory} \cite{Virtual1999} published in 1999. 
 The theme of  Turaev's paper is the search for homotopy invariants of virtual strings.  In analogy with the theory of knots, he introduced a class of slice virtual strings. Such a string underlies
a closed curve on a closed surface which is contractible in a handlebody bounded by this surface. He defined a one-variable polynomial  that encodes the obstructions for two virtual 
strings to be homotopic and for a string to be slice. He then introduced a Lie cobracket in the free abelian group generated by the homotopy classes of strings which  gives rise to
 a Hopf algebra structure in the polynomial algebra generated by the homotopy classes of strings. Using the notion of virtual string, Turaev gave a complete combinatorial 
 description of closed curves on the 2-sphere in terms
of Gauss words and bipartitions. 
 
The year after, motivated by the problem of finding a geometric presentation of words by closed curves on surfaces, Turaev, in his paper \emph{Curves on surfaces, charts, and words} \cite{Curves2005} (2005),  was led to the problem of classifying curves on surfaces. In the same paper, he gave a  combinatorial description of  curves on oriented surfaces in terms of certain
permutations called charts. He described automorphisms of curves in terms of charts, he computed the total
number of topological classes of curves counted with appropriate weights and he discussed relations between curves, words, and
complex structures on surfaces.

 In his paper \emph{Coalgebras of words and phrases}  \cite{Coalgebras2007} (2007), Turaev reformulated the theory of tensor algebras in terms of words and phrases.
 In the paper \emph{Cobordisms of words}  \cite{Cobordisms2008} (2008), he studied  an equivalence relation for words, which he called cobordism. He produced cobordism invariants of words and a measure of how far two non-cobordant nanowords are from being cobordant. In doing so, he studied surgeries of nanowords, inspired from the surgery operation of manifolds. He defined a notion of symmetry for sub-nanophrases, which is an analogue of Poincar\'e duality for manifolds.  With this relation, symmetric nanophrases become 
 the analogues of manifolds.

\subsection{Enumeration problems in topology and group theory}

  Turaev worked on other fundamental problems in classical topology, and I would like to mention here his work on the existence problem of sections of fibrations. 

In the paper \cite{enumeration2009} (2009), he gives the solution of a certain number of enumeration problems for homotopy classes of sections of locally trivial fiber bundles over surfaces. In particular, he gives a formula for the number of sections up to equivalence (homotopy and action of the second homotopy group of the fiber) of a locally trivial fiber bundle over a surface, under the hypothesis that  the fiber is connected and has finite fundamental group. The formula involves certain 2-dimensional cohomology classes associated with irreducible linear complex representations of the fundamental group of the fiber. Several applications are stated. The arguments are based on the enumeration of certain classes of group homomorphisms  and a computation of the non-abelian 1-cohomology of the fundamental group of the surface. 
There are interesting group-theoretic applications, namely, to the question of the existence of commuting elements in finite groups, and that of commuting lifts.  For any two commuting elements $a$ and $b$ in a group $G$ and for any  group epimorphism $q:G'\to G$, Turaev gives a formula for the number of commuting lifts of $a,b$ to $G'$ in terms of representations of Ker $q$. In the paper \emph{Enumeration of lifts of commuting elements of a group} \cite{Natapov}, written with M. Natapov, Turaev applies this formula to several epimorphisms $q$ with the same kernel.  
They give explicit and optimal results in the case of the quaternion group of order 8.

In the paper \cite{Sections2010} which we already mentioned, Turaev studies existence and enumeration problem
for sections of Serre fibrations over compact orientable surfaces, and he obtains a complete solution of this problem under the assumption that the
fundamental group of the fiber is finite. This solution is given again in terms
of 2-dimensional cohomology classes associated with certain irreducible representations
of this fundamental group and the proofs are based on TQFT.

  \subsection{Generalizations of the Thurston norm} In his paper \cite{Thurston}, Thurston introduced a (semi)-norm on the second homology group of a 3-manifold which turned out to be a central object in the topology of 3-manifolds. This norm
   measures the complexity of a homology class as the minimal genus of a surface representing it \cite{Thurston}. 
Motivated by this definition, Turaev introduced in his paper \emph{A norm for the cohomology of 2-complexes} \cite{Norm2002} (2002) a semi-norm on the first real cohomology group of a finite 2-complex where,  in the definition of this norm,
instead of surfaces, one takes  graphs in the 2-complex.
   He also introduced a generalization of the usual Alexander--Fox polynomial of a group by introducing a notion of twisted Alexander--Fox polynomial and he showed that they determine norms on the real 1-cohomology of groups. He proved that for any finite 2-complex, the norm defined by graphs is an upper bound for the  Alexander--Fox norms derived from the fundamental group of the complex. He highlighted a parallel between this result and  the classical Seifert inequality in knot theory which says that the genus of a knot in the 3-sphere is at least half of the span of its Alexander polynomial. 
      
   Turaev's norm may be defined for manifolds, but in general, it does not agree there with Thurston's norm. (Recall that Thurston's norm is defined using complexities
of surfaces, whereas the Turaev norm is defined using complexities of graphs.) In the paper \cite{FSW}, Friedl, Silver and Williams show that in the case where the manifold is the exterior of a link in a rational homology 3-sphere,  the Thurston norm agrees with a suitable variation of Turaev's norm defined on any 2-skeleton of this manifold.

    In another paper published in 2007 and titled \emph{A function on the homology of 3-manifolds}, again in analogy with Thurston's norm, Turaev introduced, for any orientable 3-manifold $M$, a 
(semi)-norm on the homology group $H_2 (M; \mathbb{Q}/\mathbb{Z})$ which measures the complexity of a homology class. The definition is based on the notion of complexities of ``folded" surfaces representing  the homology class, that is, surfaces whose singularities are (locally)   unions of half-planes in 3-space with a common boundary line.  

One advantage of Turaev's norm is that it gives interesting information on rational homology 3-spheres, whereas Thurston's norm does not.

\subsection{Knotoids, knots, links and braids}

In his paper  \emph{Knotoids}  \cite{Knotoids2012} (2012), Turaev introduced the notion of knotoid,  an equivalence class of generic immersions of the unit interval in a surface, together with an information on the crossing (over or under). The equivalence relation is defined by Reidemeister moves away from endpoints and by isotopy of diagrams. A knotoid is an ``open knot", and in some sense, a  rudimentary version of a knot (a knot may be obtained from a knotoid in the plane by connecting its endpoints with arcs that are positioned entirely above or below the given diagram). Turaev extended the diagrammatic theory of knots to that of knotoids.  He introduced knotoid groups and studied their properties, adapting to this setting notions like the Kauffman bracket, skein relations and skein modules and algebras.

The theory of knotoids was further extended by G\"ug\"umc\"u and Kauffman in \cite{GK}. It was applied in biology; see the paper \cite{GDBS}   by Goundaroulis, Dorier, Benedetti and Stasiak in which these authors indicate that knotoids are more natural and  more useful than knots in the production of fingerprints of protein chains.

In the paper \emph{Invariants of knots and 3-manifolds from quantum groupoids} \cite{Invariants1999} (1999), Nikshych, Turaev and Vainerman used the categories of representations of finite-dimensional quantum groupoids to construct ribbon and modular categories that give rise to invariants of knots and 3-manifolds.

In the paper \emph{Cobordism of knots on surfaces}  \cite{Cobordism2008}  (2008), Turaev studied the relation of cobordism of knots in 3-manifolds that are the product of a surface with an interval. He introduced a certain number of algebraic and geometric non-trivial cobordism invariants of these knots.

In two papers titled \emph{Higher skein modules} and \emph{Higher skein modules. II} \cite{Andersen1, Andersen2} (1999 and 2001), Turaev and J. E. Andersen, inspired by the theory of Vassiliev link invariants,  developed the theory of higher Conway skein modules and higher HOMFLY skein modules.

Turaev published two papers with D. Cimasoni  titled \emph{A Lagrangian representation of tangles} \cite{CT} and \emph{A Lagrangian representation of tangles. II} \cite{CT2} in which they generalized the Bureau representation of the braid group to the category of tangles.
In another paper titled \emph{A generalization of several classical invariants of links} \cite{Generalization2007}  (2007), the same authors studied links in quasi-cylinders. Here, a quasi-cylinder over a commutative ring $R$ is an oriented 3-manifold $M$ equipped with a submodule $V$ of the $R$-module $H_1(\partial M;R)$ such that the inclusion homomorphism $V\to H_1(M;R)$ is an isomorphism. A manifold which is the product of an oriented surface and an interval is naturally  equipped with a structure of a quasi-cylinder over $\mathbb{Z}$. In their paper, Turaev and Cimasoni extended several classical invariants of links in the 3-sphere  to this new setting (linking number, the Seifert form, the Alexander--Conway polynomial, link genus, 
the Murasugi--Tristram--Levine signature of a link in Euclidean 3-space, and other invariants).

 Turaev's book with Kassel,  \emph{Braid groups}, published in 2008 \cite{KT}, is a systematic account of braid groups, containing an exposition of several important notions related to braids that were recently developed, including Dehornoy's discovery of a natural order on braid groups made in 1991 and the linearity results for these groups established by Krammer and Bigelow. The book  became a standard reference on the topic.

\subsection{Intersection of loops on surfaces}\label{s:2-intersection}

The results of the paper \emph{Fox pairings and generalized Dehn twists} \cite{Fox2013}  (2013) by Turaev and Massuyeau are
inspired by work of Kawazumi and Kuno whose aim is to generalize, in an algebraic setting, the notion of Dehn twist to non-simple curves on surfaces. In this paper, Turaev and Massuyeau introduced a notion of Fox pairing to define automorphisms of Malcev completions of groups. (The notion of Malcev completion of a group is reviewed in the chapter by Kuno, Massuyeau and Tsuji in the present volumes \cite{KMT}, which  is a survey of generalized Dehn twists and their various applications.)
An archetypal example in this setting is the homotopy intersection form introduced by Turaev in his  1978 paper \emph{Intersections of loops in two-dimensional manifolds} \cite{Intersections1978}, of which Turaev and Massuyeau  gave a tensorial description.

 In their paper \emph{The logarithms of Dehn twists} \cite{Logarithms2014} (published in 2014 but circulated in a preprint form since 2010), Kawazumi and Kuno introduced another invariant of loops on compact oriented surfaces with one boundary component which they used in their generalization of the action of Dehn twists on the
completed group ring of the fundamental group of the surface. Their work is based on a homological interpretation of the Goldman Lie algebra in the
setting of Kontsevich's formal symplectic geometry. They also used a notion of symplectic expansion of the fundamental group due to Massuyeau (see \cite{Infinitesimal2012} published in 2012, in preprint form in 2008). The collaborations between Kawazumi and Kuno and between Massuyeau and Turaev are part of a series of mutually influenced works which we describe in the next subsection.

  \subsection{Cobrackets} \label{cobracket2}

 The 1989 paper  by Turaev,
\emph{Algebras of loops on surfaces, algebras of knots, and quantization} \cite{Algebras1989},  in which he introduced his Lie cobracket on the free module generated by the homotopy classes of loops on an oriented surface, gave rise to a series of works by him and by others, that are spread over the last 30 years. The bialgebra is called Goldman--Turaev bialgebra. (It seems that the expression was first used by Kohno, and, after that, by Kawazumi and his collaborators). This bialgebra appears in several geometric and topological contexts, including the Poisson geometry of moduli spaces and the study of quantum invariants. Kawazumi, Kuno and Tsuji used it to develop a geometric approach to the so-called  Johnson--Morita theory for the Torelli group; see the survey \cite{Intersection2016}. In the paper \emph{ Intersection of curves on surfaces and their applications to mapping class groups} \cite{Intersection2015} (2015) by Kawazumi and Kuno, the Goldman--Turaev Lie bialgebra is used in a beautiful way in the study of the higher Johnson homomorphisms.
 
Several years later, Turaev and Massuyeau,  in a joint paper titled \emph{Quasi-Poisson structures on representation spaces of surfaces} \cite{Quasi2014} (2014), showed that,  for each $N\geq 1$,
 there is a canonical Poisson-type structure on the   space of homomorphisms of the fundamental group of a surface with boundary in $\mathrm{GL}_N(\mathbb{R})$. This structure determines the Poisson structure induced by the symplectic structure on the quotient of this space by the action of $\mathrm{GL}_N(\mathbb{R})$ defined by Fock and Rosly \cite{FR}, which is an analogue for surfaces with boundary of the symplectic structure on the moduli space of  representations. The work of Fock and Rosly was aimed to extend the work of Goldman \cite{Goldman1984} (1984) that gives an interpretation of the Poisson structure induced by the Weil--Petersson synmplectic form of Teichm\"uller space. 

In a subsequent paper titled \emph{Brackets in representation algebras of Hopf algebras} \cite{Brackets2018} and published in 2018, Turaev and  Massuyeau defined a Poisson bracket on a wide class of commutative algebras, realizing an algebraic generalization of the Atiyah--Bott--Goldman Poisson structures on moduli spaces of representations of surface groups. 

 In another work  titled \emph{Brackets in the Pontryagin algebras of manifolds} \cite{Brackets2017} (2017), Massuyeau and Turaev studied the Pontryagin algebra of smooth oriented manifolds with non-empty boundary and they constructed
a bracket in the associated representation algebras. This work belongs to the field of ``string topology" inaugurated by
 Chas and Sullivan, who introduced this name after they generalized Goldman's Lie bracket to manifolds of arbitrary dimension. In the case of surfaces, the new bracket coincides with the quasi-Poisson bracket defined on the space of homomorphisms from the fundamental group of the surface to $\mathrm{GL}_N(\mathbb{R})$. In dimension $\geq 3$, the representation algebras are graded and the bracket satisfies the axioms of a Poisson bracket with appropriate signs.

In the papers \cite{Higher2017} and \cite{Goldman-Turaev2018}  by Alekseev,  Kawazumi,  Kuno and Naef, a higher-genus Kashiwara--Vergne problem is introduced using the Turaev cobracket and an isomorphism is established between the Goldman--Turaev Lie bialgebra and another Lie bialgebra structure arising from this Kashiwara--Vergne problem. In the case of surfaces of genus zero, a similar result was obtained by Massuyeau in \cite{Formal2018}, using the Kontsevich integral.  A review of these developments is contained in the chapter by Kawazumi in the present volume \cite{Kawazumi-Turaev}.

In the two  recent papers, \emph{Loops in surfaces and star-filling} \cite{Loops2019}, written in 2019, and  \emph{Quasi-Lie bialgebras of loops in quasi-surfaces} \cite{Quasi2020}, written in 2020, Turaev gives new points of view on the algebraic intersection operations of homotopy classes of loops on surfaces.  In the first paper, he uses a notion of filling a surface with boundary by a certain graph which is star-like and which is called a \emph{star}. He establishes relations between his work and that of Kawazumi--Kuno in \cite{Logarithms2014} and his work with Massuyeau in \cite{Quasi2014}. In the second paper he develops a notion he calls \emph{quasi-surface}. This is an oriented surface with boundary to which has been glued, along a set of disjoint segments on the boundary (called the ``gates"), a certain number of topological spaces. He introduces a generalization of the Lie bialgebra formed by Goldman's bracket and his cobracket to the setting of quasi-surfaces.
He introduces a structure he calls quasi-Lie bialgebra, and operations on loops on quasi-surfaces that satisfy the axioms of such a structure. Finally, in his paper \emph{Quasi-Poisson structures on moduli space of quasi-surfaces} \cite{QP2020}, he develops a generalization of the classical Atiyah--Bott Poisson bracket on the moduli spaces of surfaces to quasi-Poisson brackets on the moduli spaces of quasi-surfaces.

\subsection{Combinatorial group theory, metric geometry and phylogenetics}
 
 In 2014, Turaev published a paper titled \emph{Matching groups and gliding systems} \cite{Gliding2014}, in which he makes a relationship between what he calls \emph{gliding systems} in groups and non-positively curved cube complexes. Gliding systems are tools for constructing nonpositively curved complexes and groups. Using such systems, Turaev introduced a notion of \emph{matching groups} of which he gives an interpretation as   fundamental groups of nonpositively curved cubed complexes. He showed that these groups are torsion-free, residually nilpotent, residually finite, biorderable, biautomatic, have solvable word and conjugacy problems, and that they satisfy the Tits alternative. He also proved that they embed in $\mathrm{SL}(n,\mathbb{Z})$ for some $n$, and in finitely generated right-handed
Artin groups. 
  
On metric geometry, Turaev wrote the two papers \emph{Trimming of finite metric spaces} \cite{trim1} (2016) and \emph{Trimming of metric spaces and the tight span} \cite{trim2} (2018) 
 in which he introduced and studied the notion of  \emph{trim metric spaces} (or trim pseudo-metric spaces). This theory establishes a coherent setting for finite metric spaces, metric graphs and metric trees. It is in the tradition of the geometric theories of metric spaces developed in the 20th century by K. Menger, A.D. Alexandrov and H. Busemann.

We recall now that phylogenetics is the branch of genetics where one studies genetical modifications in living species (animal of vegetal), and in particular, the transformations and the relations among them. These relations are naturally presented in the form of trees, called phylogenetic trees, that encode the history evolution or relation and the degree of parentage among individuals.

Turaev used his theory of trimming in a paper on phylogenetics he wrote in 2018, titled \emph{Axiomatic phylogenetics} \cite{Phylo}.
In this paper, using the language of quivers, he introduced a system of axioms for a mathematical  evolution theory, with an emphasis on phylogenetic. We recall that a quiver is a directed graph where loops and multiple edges
between vertices are allowed. Turaev says that a vertex  of a quiver is phylogenetic if all possible evolutions from primitive vertices to the given vertex  have a common part. This part is viewed as the canonical evolutionary history of this vertex. He introduces then the concept of phylogenetic quiver, in which all vertices are phylogenetic and all edges are non-degenerate in a certain sense.
The principal aim of his study of phylogenetic quivers is their construction
and classification. He gives examples of phylogenetic quivers arising in various branches of mathematics: set theory,
group theory, and the theory of metric spaces.  In the last section of his paper, he indicates how trim metric spaces  can be used in phylogenetics.

\section{In lieu of a conclusion}

In this short survey, I tried to convey the fact that Turaev's ideas are now at the cutting edge of current research. A certain number of dedicated mathematicians are working on problems he formulated, using notions he discovered,  in a field sometimes called ``quantum topology", whose central core is topology, enhanced by algebra and inspired by theoretical physics.

Turaev's style in mathematics is elegant, characterized by its clarity, the emphasis on essential ideas, and at the same time by a meticulousness for the details. In his writings, he establishes connections between beautiful concepts and, whenever this is useful, he likes to add historical comments. 

His style in life is essentially the same: clarity, faithfulness and going to the essential. He is knowledgeable in literature, art and history, and extremely open to new ideas. He arrived to France without any knowledge of French, and his French became perfect in an amazingly short period of time.  I remember that when he decided to emigrate to France, he was anxious only about two things:  first, moving his parents and his son from the Soviet Union to Western Europe, and, then, sending his books from Leningrad to Strasbourg. He had a durable influence, sometimes with a profound effect, on  people who were close to him. He strongly values the notion of friendship. I learned, in part through my relation with him, that mathematics is also a question of friendship.

\bigskip

\noindent{\it Acknowledgements.} I would like to thank Christian Blanchet, Nariya Kawazumi and Gwenaël Massuyeau who read a preliminary version of this article.


\begin{thebibliography}{99}
   
   
\bibitem{Hyperbolicity}  C. Adams, O. Eisenberg, J. Greenberg, K. Kapoor, Z. Liang,  K. O'Connor,  N. Pacheco-Tallaj and Y. Wang, Turaev hyperbolicity of classical and virtual knots, preprint, 2019.
 

  
\bibitem{Higher2017} A. Alekseev,  N. Kawazumi,  Y. Kuno and F. Naef,   Higher genus Kashiwara--Vergne problems and the Goldman-Turaev Lie bialgebra.  
C. R. Math. Acad. Sci. Paris 355 (2017), no. 2, 123127.
 
\bibitem{Goldman-Turaev2018} A. Alekseev,  N. Kawazumi,  Y. Kuno and F. Naef,  
The Goldman-Turaev Lie bialgebra in genus zero and the Kashiwara--Vergne problem.  Adv. Math. 326  (2018), 1-53. 



\bibitem{Andersen1} J. E. Andersen and V. G.  Turaev,  Andersen, 
Higher skein modules. 
J. Knot Theory Ramifications 8 (1999), no. 8, 963-984.
\bibitem{Andersen2} J. E. Andersen and V. G.  Turaev,  
Higher skein modules. II. In: Topology, ergodic theory, real algebraic geometry, 21-30. In;
Topology, ergodic theory, real algebraic geometry.
Rokhlin's memorial. Edited by V. Turaev and A. Vershik. American Mathematical Society Translations, Series 2, 202. Advances in the Mathematical Sciences, 50. American Mathematical Society, Providence, RI, 2001. 


 
\bibitem{Atiyah-Work} M. F. Atiyah, On the work of Edward Witten, 
Proc. Int. Congr. Math., Kyoto/Japan 1990,  (1991) Vol. I, 31-35.


\bibitem{BDR} C. Blanchet and M. De Renzi,  Modular categories and TQFTs beyond semisimplicity, this volume, p. ???


\bibitem{Brody} E. J. Brody, The topological classification of the lens spaces, Annals of Mathematics,
Second Series, 71 (1960) No. 1, 163-184 



\bibitem{CT} D. Cimasoni, and V. G. Turaev,  A Lagrangian representation of tangles. Topology 44 (2005), no. 4, 747-767.  

\bibitem{CT2} D. Cimasoni, and V. G. Turaev,   A Lagrangian representation of tangles. II. Fund. Math. 190 (2006), 11-27. 


\bibitem{Generalization2007} D. Cimasoni and V. G. Turaev,  
A generalization of several classical invariants of links.  
Osaka J. Math. 44 (2007), no. 3, 531-561.


\bibitem{Connes-Kreimer1998} A. Connes and D. Kreimer, Hopf algebras, renormalization and noncommutative geometry. Comm. Math. Phys. 199 (1998), no. 1, 203-242.  
 
 \bibitem{CoT} F. Costantino and D. Thurston, 
3-manifolds efficiently bound 4-manifolds. J. Topol. 1 (2008), no. 3, 703-745.
 
\bibitem{DL} O. T. Dasbach and A. M. Lowrance,    A Turaev surface approach to Khovanov homology, Quantum Topol. 5 (2014), 425-486.

 \bibitem{Drinfeld1983} V. G. Drinfeld,  Hamiltonian structures on Lie groups, Lie bialgebras and the geometric meaning of classical Yang--Baxter equations.  Dokl. Akad. Nauk SSSR 268 (1983), no. 2, 285-287.
    
    


 \bibitem{DW} R. H. Dijkgraaf and E. Witten, Topological gauge theories and group cohomology. Comm. Math. Phys. 129 (1990), no. 2, 393-429.
 
        \bibitem{EN} H. Endo and S. Nagami, Signature of relations in mapping class groups and nonholomorphic
Lefschetz fibrations, Trans. Amer. Math. Soc. 357 (2005), no. 8, 3179-3199.
        
\bibitem{Absolute1999} M. Farber and V. G. Turaev, 
Absolute torsion. Tel Aviv Topology Conference: Rothenberg Festschrift (1998), 73-85,
Contemp. Math., 231, Amer. Math. Soc., Providence, RI, 1999.
 


 \bibitem{Poincare2000} M. Farber and V. G. Turaev, Poincar\'e--Reidemeister metric, Euler structures, and torsion. J. Reine Angew. Math. 520 (2000), 195-225.

\bibitem{FR} V. V. Fock, and A. A. Rosly,  
Poisson structure on moduli of flat connections on Riemann surfaces and the r-matrix.   Moscow Seminar in Mathematical Physics, 67-86,
Amer. Math. Soc. Transl. Ser. 2, 191, Adv. Math. Sci., 43, Amer. Math. Soc., Providence, RI, 1999. 
 
\bibitem{T-Vesnin} E. A. Fominykh, V. G. Turaev  and  A. Yu. Vesnin, Three-dimensional manifolds with poor spines, Proceedings of the Steklov Institute of Mathematics, 288 (2015),  29-38


\bibitem{Franz} W. Franz, \"Uber die Torsion einer \"Uberdeckung, J. Reine Angew. Math. 173 (1935), 245-254.


      \bibitem{FSW} S. Friedl, D. S. Silver and S. G. Williams,
  The Turaev and Thurston norms, Pacific J. Math. 284 (2016), no. 2, 365-382. 


\bibitem{Tetrahedral2012} N. Geer,  R. Kashaev and V. G. Turaev,  Tetrahedral forms in monoidal categories and 3-manifold invariants.  
J. Reine Angew. Math. 673 (2012), 69-123.
 


   
 \bibitem{Modified2009} N. Geer,  B. Patureau-Mirand and V. G. Turaev,  
Modified quantum dimensions and re-normalized link invariants. 
Compos. Math. 145 (2009), no. 1, 196-212.


\bibitem{Modified2011} N. Geer,  B. Patureau-Mirand and V. G. Turaev,  
Modified $6j$-symbols and 3-manifold invariants.  
Adv. Math. 228 (2011), no. 2, 1163-1202.
 
     
\bibitem{Gauss} C. F. Gauss, Werke, Vol. VIII, Teubner, Leipzig, 1900.




\bibitem{Goldman1984}  W. M. Goldman,  The symplectic nature of fundamental groups of surfaces.
Adv. in Math. 54 (1984), no. 2, 200-225. 
  
  

 \bibitem{GDBS}    D. Goundaroulis, J. Dorier, F. Benedetti and A. Stasiak, Studies of global and local entanglements of individual protein chains using the concept of
knotoids, Scientific Reports (electronic journal) vol. 7,  (2017), Article number,  6309


\bibitem{GK} N. G\"ug\"umc\"u, L. H. Kauffman, New invariants of knotoids, European Journal of Combinatorics 65 (2017),  186-229. 


\bibitem{KT} C. Kassel and V. G. Turaev,  
Braid groups.
  Graduate Texts in Mathematics, 247. Springer, New York, 2008. 
  



\bibitem{State1987} L. H. Kauffman, State models and the Jones polynomial. Topology 26 (1987), no. 3, 395-407. 

\bibitem{Virtual1999} L. H. Kauffman, Virtual knot theory. European J. Combin. 20 (1999), no. 7, 663-690.

\bibitem{Kawazumi-Turaev} N. Kawazumi, Some algebraic aspects of the Turaev cobracket, this volume, p. ???


 \bibitem{Intersection2015}  N. Kawazumi and Y. Kuno, Intersection of curves on surfaces and their applications to mapping class groups. Ann. Inst. Fourier (Grenoble) 65 (2015), no. 6, 2711-2762.


 \bibitem{Intersection2016}  N. Kawazumi and Y. Kuno,
The Goldman-Turaev Lie bialgebra and the Johnson homomorphisms. Handbook of Teichm\"uller theory (ed. A. Papadopoulos). Vol. V, 97-165,
IRMA Lect. Math. Theor. Phys., 26, Eur. Math. Soc., Z\"urich, 2016. 


 \bibitem{Logarithms2014}  N. Kawazumi and Y. Kuno,  
The logarithms of Dehn twists. 
Quantum Topol. 5 (2014), no. 3, 347-423.
 
 


\bibitem{KK} S. Kim and I. Kofman, Turaev Surfaces, In: A Concise Encyclopedia of Knot Theory,  C. Adams, E. Flapan, A. Henrich, L. Kauffman, L. Ludwig, and S. Nelson (ed.)
 Taylor \& Francis, to appear.

 \bibitem{King}  S. King, S. V. Matveev,  V. Tarkaev, V. G. Turaev, Dijkgraaf-Witten  $\mathbb{Z}/2$-invariants for Seifert manifolds. J. Knot Theory Ramifications 26 (2017), no. 1,  7 pages.

    \bibitem{Kirby} R. Kirby, The calculus of framed links in $S^3$. Invent. Math. 45 (1978), 35-56.
    
\bibitem{KMT} Y. Kuno, G. Massuyeau and S. Tsuji, Generalized Dehn twists in low-dimensional topology, this volume, p. ???


 \bibitem{Kuperberg} G. Kuperberg,  Quantum invariants of knots and 3-manifolds, Bulletin of the AMS, Volume 33, Number 1, January 1996, 107-110.
        


\bibitem{Infinitesimal2012} G. Massuyeau, 
Infinitesimal Morita homomorphisms and the tree-level of the LMO invariant. 
Bull. Soc. Math. France 140 (2012), no. 1, 101-161. 


\bibitem{Formal2018} G. Massuyeau, 
Formal descriptions of Turaev's loop operations.  
Quantum Topol. 9 (2018), no. 1, 39-117.



\bibitem{Fox2013} G. Massuyeau and V. G. Turaev, 
Fox pairings and generalized Dehn twists.  Ann. Inst. Fourier 63  (2013), No. 6, 2403-2456. 
 
 
\bibitem{Quasi2014} G. Massuyeau and V. G. Turaev,  Quasi-Poisson structures on representation spaces of surfaces. Int. Math. Res. Not. IMRN 2014, no. 1, 1-64.


\bibitem{Brackets2017} G. Massuyeau and V. G. Turaev, Brackets in the Pontryagin algebras of manifolds. M\'em. Soc. Math. Fr. (N.S.) No. 154, 2017, 138 pp.
  
  
\bibitem{Brackets2018} G. Massuyeau and V. G. Turaev,  Brackets in representation algebras of Hopf algebras. J. Noncommut. Geom. 12 (2018), no. 2, 577-636.


\bibitem{Matveev} S. V.  Matveev, Algorithmic Topology and Classification of 3-Manifolds, Springer, Berlin,
2007.


\bibitem{Matveev1} S. V. Matveev, V. G. Turaev,  Dijkgraaf-Witten invariants over $\mathbb{Z}/2$ for 3-manifolds.  
Dokl. Akad. Nauk 460 (2015), no. 1, 15-17. English translation in
Dokl. Math. 91 (2015), no. 1, 9-11
 


\bibitem{Milnor1957} J. Milnor,  Isotopy of links. Algebraic geometry and topology. A symposium in honor of S. Lefschetz, p. 280-306. Princeton University Press, Princeton, N. J., 1957.

\bibitem{Milnor1962} J. Milnor, A duality theorem for Reidemeister torsion.
Ann. of Math. (2) 76 (1962), 137-147. 


 
\bibitem{Natapov} M. Natapov and V. G. Turaev,  
Enumeration of lifts of commuting elements of a group. 
J. Algebra 324 (2010), no. 10, 2732-2741.

 


\bibitem{Invariants1999} D. Nikshych, V. G. Turaev and L. Vainerman, Invariants of knots and 3-manifolds from quantum groupoids.  
Proceedings of the Pacific Institute for the Mathematical Sciences Workshop "Invariants of Three-Manifolds'' (Calgary, AB, 1999).
Topology Appl. 127 (2003), no. 1-2, 91-123.
 

 
\bibitem{Planar1975}  C. D. Papakyriakopoulos, Planar regular coverings of orientable closed surfaces. Knots, groups, and 3-manifolds (Papers dedicated to the memory of R. H. Fox), 261-292. Ann. of Math. Studies, No. 84, Princeton Univ. Press, Princeton, N.J., 1975.  

 
  \bibitem{Porter2007} T. Porter,  Formal homotopy quantum field theories. II. Simplicial formal maps. Categories in algebra, geometry and mathematical physics, 375-403, Contemp. Math., 431, Amer. Math. Soc., Providence, RI, 2007.
  
  
\bibitem{Porter2012} T. Porter, Homotopy quantum field theory (Book Review), Bull. Amer. Math. Soc. (N.S.) 49 (2012), no. 2, 337-345.  
 


   \bibitem{PT2007} T. Porter and V. G. Turaev,  
Formal homotopy quantum field theories. I. Formal maps and crossed $\mathcal{C}$-algebras. 
J. Homotopy Relat. Struct. 3 (2008), no. 1, 113-159. 
 



\bibitem{Przytycki1991} J. Przytycki, Skein modules of 3-manifolds, Bull. Polish Acad. Sci. Math. 39 (1991), 91-100.
100.


\bibitem{R1} K. Reidemeister, Die Klassifikation der Linsenräume, Abh. Math. Sem. Hamburg 11 (1935), 102-109.


\bibitem{R1988} N. Yu. Reshetikhin, Quantized universal enveloping algebras, the Yang--Baxter equation, and invariants of links', Akad. Nauk SSSR Leningrad. Otdel. Mat. Inst. Steklov. Preprint No. E-4-87, 1988.

  \bibitem{Ribbon1990} N. Yu. Reshetikhin and V. G. Turaev, 
Ribbon graphs and their invariants derived from quantum groups.  
Commun. Math. Phys. 127  (1990), No. 1, 1-26. 


 \bibitem{Invariants1991} N. Yu. Reshetikhin and V. G. Turaev, 
Invariants of 3-manifolds via link polynomials and quantum groups.  
Invent. Math. 103  (1991), No. 3, 547-597.       



 \bibitem{ST2010}  M. D. Staic and V. G. Turaev,  Remarks on 2-dimensional HQFTs. Algebr. Geom. Topol. 10 (2010), no. 3, 1367-1393. 
 

      
      
\bibitem{Thurston} W. Thurston, A norm for the homology of 3-manifolds. Mem. Amer. Math. Soc. 59 (1986).


   \bibitem{Reidemeister1976} V. G. Turaev, Reidemeister torsion and the
Alexander polynomial, Mat. Sb. (N.S.)
Volume 101(143) (1976), Number 2(10), 252-270. English translation Math. USSR Sb. 30 (1976) No. 2, 220-237.
      
\bibitem{Milnor1976} V. G. Turaev, Milnor invariants and Massey products.  
Studies in topology, II. Zap. Nauchn. Sem. LOMI, 66 (1976), 189-203. English transl. J. Soviet Math. 12 (1979), 128-137.
 
 \bibitem{Intersections1978} V. G. Turaev, Intersections of loops in two-dimensional manifolds.  
 Mat. Sb. 106 (148) (1978), no. 4, 566-588. English translation Math. USSR, Sb. 35, 229-250 (1979). 
 
 \bibitem{Cohomology1983} V. G. Turaev, Cohomology rings, linking coefficient forms and invariants of spin structures in three-dimensional manifolds. Mat. Sb. (N.S.) 120(162) (1983), no. 1, 68-83.
 
  \bibitem{Cocycle1984} V. G. Turaev, A cocycle of the symplectic first Chern class and Maslov indices,  Funktsional. Anal. i Prilozhen. 18 (1984), no. 1, 43-48. English translation in Functional analysis and its appl. 118 (1984), 35-39.
  
  
    \bibitem{Reidemeister1986} V. G. Turaev, Reidemeister torsion in knot theory. Russian Math. Surveys 41 (1986) 119-182.
 
  \bibitem{First1987} V. G. Turaev, The first symplectic Chern class and Maslov indices. Studies in topology, V. Zap. Nauchn. Sem. Leningrad. Otdel. Mat. Inst. Steklov. (LOMI) 143 (1985), 110-129, 178. English translation in J. of Soviet Math. 37 (1987), 1115-1127.
  
  

\bibitem{Simple1987} V. G. Turaev. A simple proof of the Murasugi and Kauffman theorems on alternating links. Enseign. Math. 33 (1987), 203-225.


\bibitem{Turaev1988} V. G. Turaev, The Conway and Kauffman modules of a solid torus (Russian). Notes of LOMI
scientific seminars 167 (1988), 79-89. English translation: J. Soviet Math. 52 (1990), 2799-2805.

        
  \bibitem{Yang1988} V. G. Turaev, 
The Yang--Baxter equation and invariants of links. 
Invent. Math. 92  (1988), No. 3, 527-553.      

  \bibitem{Classification1988} V. G. Turaev,  Classification of oriented Montesinos links via spin structures. Topology and geometry: Rokhlin Seminar, 271-289, Lecture Notes in Math., 1346, Springer, Berlin, 1988.
      
      

\bibitem{Towards1988} V. G. Turaev, 
Towards the topological classification of geometric 3-manifolds. Topology and geometry, Rokhlin Seminar, 291-323,
Lecture Notes in Math., 1346, Springer, Berlin, 1988. 



\bibitem{Algebras1989}  V. G. Turaev,
Algebras of loops on surfaces, algebras of knots, and quantization. In: Braid group, knot theory and statistical mechanics, 59-95,
Adv. Ser. Math. Phys., 9, World Sci. Publ., Teaneck, NJ, 1989.



\bibitem{Euler1989} V. G. Turaev, Euler structures, nonsingular vector fields, and Reidemeister-type torsions.  
Izv. Akad. Nauk SSSR Ser. Mat. 53 (1989), no. 3, 607-643.



\bibitem{Three-dimensional1989}  V. G. Turaev, 
Three-dimensional Poincar\'e complexes: Homotopy classification and splitting. Mat. Sb.  (1989) 180, No. 6, 809-830. English transl. Math. USSR, Sb. 67  (1990), No. 1, 261-282.


\bibitem{Operator invariants1989} V. G. Turaev, Operator invariants of tangles and $R$-matrices, Izv. Akad. Nauk SSSR Ser. Mat. (1989), Vol. 53,	Issue 5, 1073-1107. English translation:  Mathematics of the USSR-Izvestiya (1990) 35:2, 411-444.



   \bibitem{State1991} V. G. Turaev, 
State sum models in low-dimensional topology. 
Proc. Int. Congr. Math., Kyoto, 1990,  (1991) Vol. I, 689-698. 


\bibitem{Skein1991}  V. G. Turaev, 
Skein quantization of Poisson algebras of loops on surfaces.  
Ann. Sci. \'Ec. Norm. Sup\'er. (4) 24  (1991), No. 6, 635-704. 


\bibitem{Quantum1991} V. G. Turaev, Quantum invariants of 3-manifolds and a glimpse of shadow topology.  
C. R. Acad. Sci., Paris, S\'er. I 313  (1991), No. 6, 395-398. 


      
\bibitem{Shadow1992}  V. G. Turaev, 
Shadow links and face models of statistical mechanics.
J. Differ. Geom. 36  (1992), No. 1, 35-74. 
 
 
\bibitem{Modular1992} V. G. Turaev, Modular categories and 3-manifold invariants.  In: 
Topological and quantum group methods in field theory and condensed matter physics.
Internat. J. Modern Phys. B 6 (1992), no. 11-12, 1807-1824.
 
 
 

\bibitem{Quantum1993-Wenzl} V. G. Turaev, and H. Wenzl,  
Quantum invariants of 3-manifolds associated with classical simple Lie algebras.
Internat. J. Math. 4 (1993), no. 2, 323-358.
  

\bibitem{Quantum1993} V. G. Turaev, Quantum invariants of links and 3-valent graphs in 3-manifolds.
Inst. Hautes Études Sci. Publ. Math. No. 77 (1993), 121-171.


 \bibitem{Axioms1994} V. G.  Turaev,  
Axioms for topological quantum field theories.  
Ann. Fac. Sci. Toulouse Math. (6) 3 (1994), no. 1, 135-152.
  
 
\bibitem{Quantum1994} V. G. Turaev, Quantum invariants of knots and 3-manifolds. 
De Gruyter Studies in Mathematics, 18. Walter de Gruyter \& Co., Berlin, 1994.  3d ed. 2016.
 
 \bibitem{Homotopy1999} V. G. Turaev, Homotopy field theory in dimension 2 and group-algebras,  	arXiv:math/9910010, 1999.
  
\bibitem{Homotopy2000} V. G. Turaev,  Homotopy field theory in dimension 3 and crossed group-categories,  	arXiv:math/0005291, 2000




\bibitem{Introduction2001} V. G. Turaev, Introduction to combinatorial torsions.  
Notes taken by Felix Schlenk. Lectures in Mathematics ETH Z\"urich. Birkhäuser Verlag, Basel, 2001. 


  
       
  
  \bibitem{Norm2002} V. G. Turaev, A norm for the cohomology of 2-complexes, Alg. Geom. Top. 2 (2002), 137-155.


\bibitem{Virtual2004} V. G. Turaev,  
Virtual strings.
Ann. Inst. Fourier 54, No. 7, 2455-2525 (2004). 
 


\bibitem{Loops2005} V. G. Turaev,
Loops on surfaces, Feynman diagrams, and trees.  
J. Geom. Phys. 53  (2005), No. 4, 461-482. 


 
\bibitem{Curves2005}   V. G. Turaev,  
Curves on surfaces, charts, and words. 
Geom. Dedicata 116 (2005), 203-236.



 \bibitem{Knots2006}  V. G. Turaev,  
Knots and words.
Int. Math. Res. Not. 2006, Art. ID 84098, 23 pp.


  \bibitem{Topology2007}   V. G. Turaev,  
Topology of words. 
Proc. Lond. Math. Soc. (3) 95 (2007), no. 2, 360-412.


 
       
        
 \bibitem{Function2007} V. G. Turaev, A function on the homology of 3-manifolds,   Algebr. Geom. Topol.
    Volume 7, Number 1 (2007), 135-156.



\bibitem{Dijkgraaf2007} V. G. Turaev,  
Dijkgraaf-Witten invariants of surfaces and projective representations of groups.  J. Geom. Phys. 57 (2007), no. 11, 2419-2430.
 




 \bibitem{Coalgebras2007}    V. G. Turaev,  
Coalgebras of words and phrases. 
J. Algebra 314 (2007), no. 1, 303-323.
 




\bibitem{Lectures2007}   V. G. Turaev,  
Lectures on topology of words. 
Jpn. J. Math. 2 (2007), no. 1, 1-39.


\bibitem{Cobordism2008} V. G. Turaev, Cobordism of knots on surfaces. 
J. Topol. 1 (2008), no. 2, 285-305.
 
  \bibitem{Cobordisms2008}  V. G. Turaev,   Cobordisms of words.  
Commun. Contemp. Math. 10 (2008), suppl. 1, 927-972.




 
\bibitem{enumeration2009} V. G. Turaev, On certain enumeration problems in two-dimensional topology.  Math. Res. Lett. 16 (2009), no. 3, 515-529.

 
\bibitem{Homotopy2010} V. G. Turaev, Homotopy quantum field theory, EMS Tracts Math., 10, Eur. Math. Soc., Z\"urich, 2010


\bibitem{Sections2010} V. G. Turaev,  
Sections of fiber bundles over surfaces and TQFTs.  
Quantum Topol. 1 (2010), no. 3, 275-319.



\bibitem{Knotoids2012} V. G. Turaev, Knotoids.  
Osaka J. Math. 49 (2012), no. 1, 195-223.


 \bibitem{Gliding2014} V. G. Turaev,  
Matching groups and gliding systems,
J. Geom. Phys. 81 (2014), 128-144. 
 



\bibitem{trim1} V. G. Turaev, Trimming of finite metric spaces, arXiv:1612.06760, 2016.
 
 

 
\bibitem{Phylo} V. G. Turaev, Axiomatic phylogenetics, arXiv:1702.06388,  2017.
  
  \bibitem{trim2}  V. G. Turaev, Trimming of metric spaces and the tight span, Discrete Mathematics 341
(2018), 2925-2937.


\bibitem{Loops2019} V. G. Turaev, Loops in surfaces and star-filling, preprint, 2019.


\bibitem{Quasi2020} V. G. Turaev,  Quasi-Lie bialgebras of loops in quasi-surfaces, preprint, 2020.
 


\bibitem{QP2020} V. G. Turaev, Quasi-Poisson structures on moduli space of quasi-surfaces, preprint, 2020
  


\bibitem{Unoriented2006} V. G. Turaev, and P.  Turner, 
Unoriented topological quantum field theory and link homology.  
Algebr. Geom. Topol. 6 (2006), 1069-1093. 


 
\bibitem{Tambara2012} V. G. Turaev, and L. Vainerman,  
The Tambara-Yamagami categories and 3-manifold invariants.  
Enseign. Math. (2) 58 (2012), no. 1-2, 131-146.


 
\bibitem{TV1} V. G. Turaev, and A. Virelizier,
On 3-dimensional homotopy quantum field theory, I. 
Internat. J. Math. 23 (2012), no. 9, 28 pp.
 


\bibitem{TV3} V. G. Turaev, and A. Virelizier, Surgery HQFT, arXiv:1303.1331, 2013


\bibitem{TV2} V. G. Turaev, and A. Virelizier,
On 3-dimensional homotopy quantum field theory, II: The surgery approach.  
Internat. J. Math. 25 (2014), no. 4, 66 pp. 

\bibitem{TV} V. G. Turaev, and A. Virelizier,   Monoidal categories and topological field theory. Progress in Mathematics, 322. Birkh\"auser/Springer, Cham, 2017. 



\bibitem{TV4} V. G. Turaev, and A. Virelizier, On 3-dimensional homotopy quantum field theory III: comparison of two approaches,  	arXiv:1911.10257, 2017.


 \bibitem{Intersections1983} V. G. Turaev and O. Ya. Viro,  
Intersections of loops in two-dimensional manifolds. II: Free loops.  
Mat. Sb. (N. S.) 121(163), No. 3, 359-369 (1983). English translation Math. USSR, Sb. 49, 357-366 (1984). 
 


         \bibitem{State1992}     V. G. Turaev and O. Ya. Viro, 
State sum invariants of 3-manifolds and quantum $6j$-symbols.  
Topology 31, No. 4, 865-902 (1992). 
\bibitem{Wall1965} C. T. C. Wall,   Finiteness conditions for CW-complexes. Ann. of Math. (2) 81 (1965), 56-69. 
     
     
     \bibitem{Walker} K. Walker, On Witten's 3-manifold Invariants, unpublished preprint, 1991.

 \bibitem{Witten1989}   E. Witten,  Quantum field theory and the Jones polynomial
Communications in Mathematical Physics,  121 (1989), 351-399.   
   
     \end{thebibliography}
\end{document}